\documentclass[12pt]{amsart} % uses LaTeX2e
\chardef\bslash=`\\ % p. 424, TeXbook
\def\pa{\makebox[.5em][l]
            {a\hspace{-.23em}\begin{picture}(0,0)
                             \put(0,0){\raisebox{-1.4ex}{`}}
                             \end{picture}
            }}%
\makeatletter
\def\verbatim{\interlinepenalty\@M \@verbatim
  \leftskip\@totalleftmargin\advance\leftskip2pc
  \frenchspacing\@vobeyspaces \@xverbatim}
\makeatother
\newtheorem{theorem}{Theorem}[section]
\newtheorem{corollary}[theorem]{Corollary}
\newtheorem{lemma}[theorem]{Lemma}
\newtheorem{proposition}[theorem]{Proposition}

\newtheorem{remark}[theorem]{Remark}

\theoremstyle{definition}

\numberwithin{equation}{section}

\newcounter{picture}

\newcommand{\CC}{{\mathbb C}}
\newcommand{\FF}{{\mathbb F}}
\newcommand{\HH}{{\mathbb H}}

\newcommand{\PP}{{\mathbb P}}

\newcommand{\RR}{{\mathbb R}}

\begin{document}

\title[]{ Crofton formulae and geodesic distance in hyperbolic spaces }

\date{July 28, 1997}
\author[G. Robertson]{Guyan Robertson}
\address{Department of Mathematics  \\
        University of Newcastle\\  NSW  2308\\ AUSTRALIA}
\email{guyan@maths.newcastle.edu.au }

\thanks{This research was supported by the Australian Research Council.} 
\thanks{ \hfill Typeset by  \AmS-\LaTeX}

\begin{abstract}
The geodesic distance between points in real hyperbolic space is a hypermetric,
and hence is a kernel negative type. The proof given here uses
an integral formula for geodesic distance, in terms of 
a measure on the space of hyperplanes.
An analogous integral formula, 
involving the space of horospheres, is given for complex 
hyperbolic space.
By contrast geodesic distance in a projective space is not of
negative type.
\end{abstract}      

\maketitle
\section*{Introduction}

Motivated by problems in harmonic analysis on homogeneous spaces,
J. Faraut and K. Harzallah showed in
\cite {far} that the geodesic distance $d$ on a real
 or complex hyperbolic space of dimension $n \ge
2$ is a kernel of negative type.   Equivalently, $\sqrt {d}$ is a Hilbert space distance [HV:
Chapter 5].  On the other hand, in the real hyperbolic plane, there is an
explicit Crofton integral formula for the length $L$ of a rectifiable curve $C$
in terms of a measure $\mu$ on the space of geodesics \cite[Section 3]{s1}, 
\cite[Section 4.4]{gel}.  The formula has the form $\int n(\gamma) d \mu (\gamma) = 2L$,
where $n (\gamma)$ is the number of times the geodesic $\gamma$ meets $C$ and
the integral is over the space of all geodesics.  (M.W. Crofton proved the
corresponding formula for the euclidean plane in 1868.)

This paper proves a Crofton formula for geodesic distances in real hyperbolic
space of dimension $n$, from which one obtains an explicit geometric proof of
the result of Faraut and Harzallah, in the real case.  In this case the measure
$\mu$ is defined on the space of totally geodesic submanifolds of codimension
one.  There is also a Crofton formula in the complex case.  However, the
measure is defined on the space of horospheres, and it is not clear how to
deduce that $d$ is of negative type.

Hyperbolic spaces $H_{\FF}^n$ can be defined over ${\RR}$, ${\CC}$, ${\HH}$ and (in
the case $n = 2$) the octonions ${\mathbb O}$.  In each case the space is two
point homogeneous:  if $d(x, y) = d(x^\prime, y^\prime)$ then there is an
isometry $g$ of $H_{\FF}^n$ such that $gx = x^\prime$, $gy = y^\prime$
\cite[Chapter IX, Proposition 5.1, p.355]{h1}.
  It turns out that this
property, together with the existence of a certain invariant measure is at the heart of
the Crofton formulae.  We note that the noncompact riemannian manifolds which
are two point homogeneous are exactly the euclidean spaces and hyperbolic spaces 
\cite{h1}.

The Crofton formula that we prove for $H_{\RR}^n$ has the consequence that
the geodesic distance $d$ is a {\em hypermetric} in the sense of \cite{a}.
It follows that any finite subset of
$H_{{\RR}}^n$, endowed with the metric $\sqrt d$, 
embeds isometrically in a euclidean sphere (Corollary \ref{spherical}).

Finally, we complete the investigation of Faraut and Harzallah \cite{far}
for symmetric spaces of rank one,
by showing explicitly that the geodesic distance on a projective space
is not of negative type, although it can be expressed by means of a Crofton
formula.

\bigskip
\section{motivation and preliminaries}

In the course of the proof of \cite [Proposition 1.4]{rs}, a geometric
argument was given to show that the euclidean distance $d(x,y)$
between points $x,y \in {\mathbb R^n}$ is a kernel of negative type.
In retrospect, the basis of that argument is seen to be the classical Crofton
formula asserting that $d(x,y)$ equals the measure of
the set of euclidean hyperplanes which meet the line segment
$[xy]$. The relevant measure is the appropriately normalized measure
on the space of hyperplanes which is invariant under isometries of
${\mathbb R^n}$. This measure lifts
to a measure $\mu$ on the space of half spaces of ${\mathbb R^n}$.
Then $d(x,y) = \mu(S_x \triangle S_y)$, where $S_x$ denotes the
set of half spaces which contain $x$. It is natural to try to extend
this argument to the rank one symmetric spaces considered by
J. Faraut and K. Harzallah \cite {far}. It turns out that
the only such spaces for which the geodesic distance is a
kernel of negative type are real or complex hyperbolic spaces,
and euclidean  spheres.    

Let $\FF$  be one of the (skew-) fields ${\RR}$, ${\CC}$ or ${\HH}$. 
Regard ${\FF}^{n+1}$ as a right vector space over $\FF$ and define 
a hermitian form on ${\FF}^{n+1}$ by means of the formula

$$\langle z, w \rangle = -
\overline {z^0} w^0 + \overline {z^1} w^1 + \cdots + \overline {z^n} w^n.$$

The hyperbolic space $H_{\FF}^n$ is the image in the
projective space $P_{\FF}^n$ of the set of negative vectors $\left \{ z \in
{\FF}^{n+1} : \langle z, z \rangle < 0 \right \}$.

It is convenient to
use the same notation for a negative vector $x \in {\FF}^{n+1}$ and its
equivalence class in $P_{\FF}^n$.  The geodesic distance between points $x,
y$ in  $X$ is given by $\cosh d(x, y) = \frac{\left |
\langle x, y \rangle \right |}{ \left ( \langle x, x \rangle \langle y, y
\rangle \right )^{\frac{1}{2}}}$.

The octionic hyperbolic plane $H_{\mathbb O}^2$ requires a more
involved definition \cite{mos}.

\bigskip
\section{real hyperbolic space}\label{rhyperbolic}

Denote by $X$ the hyperbolic space $H_{\RR}^n$ [CG].   Then $X$ is the image in the
projective space ${\RR}{\PP}^n$ of the set of negative vectors 
$$\left \{ x \in
{\RR}^{n+1} : \langle x, x \rangle = - (x^0)^2
+ (x^1)^2 + \ldots (x^n)^2 < 0 \right \}.$$    The orthogonal group $G = O(1, n)$ is the
subgroup of $GL(n+1, {\RR})$ which preserves the form $\langle \cdot , \cdot
\rangle$, and $G$ acts isometrically and transitively on $X$.  The stabilizer
in $G$ of the point $x_0 = (1, 0, \ldots , 0) \in X$ is the compact group $K =
O(1) \times O(n)$.  Thus $X$ is isomorphic to the topological homogeneous
space $G/K$.  The hyperbolic space  $H_{{\RR}}^{n-1}$ embeds naturally into $X
= H_{{\RR}}^n$ as the subspace $S_0$ consisting of all points with last
coordinate equal to zero.

Every totally geodesic submanifold of codimension one in $X$ is a
$G$--translate of $S_0$ \cite[Proposition 2.5.1]{chen}.  It is convenient simply
to refer to such a submanifold as a hyperplane.  The space ${\mathfrak S}$ of all
hyperplanes may therefore be identified with the topological homogeneous space
$G/G(S_0)$, where $G(S_0)$ is the subgroup of $G$ consisting of elements which
leave $S_0$ globally invariant.  By  \cite[Lemma 4.2.1]{chen}, we have $G(S_0)
\cong O(1, n-1) \times O(1)$.  The groups $G$ and $G(S_0)$ are both
unimodular \cite[Chapter X, Proposition 1.4]{h1}.  (A direct proof is given in
\cite[Proposition C.4.11]{bp}.)  It follows \cite[Chapter 3, p.140, Corollary
4]{nac} that there is a nonzero positive $G$-invariant measure $\mu_{\mathfrak S}$
on the space ${\mathfrak S}$ of hyperplanes.

We also consider the space ${\mathfrak H}$ of half spaces in $X$.  These are the
$G$--translates of the half--space $H_0$ consisting of points with last
coordinate positive.  The group $G(H_0) \cong O(1, n-1)$ is unimodular and there
is a corresponding invariant measure $\mu_{\mathfrak H}$ on ${\mathfrak H} = G/G(H_0)$.

There is a natural double covering $\pi : {\mathfrak H} \rightarrow {\mathfrak S}$ and so
by uniqueness of the measures (up to a positive multiple) \cite[Chapter 2, p.95,
Corollary]{nac} we may assume that $\mu_{\mathfrak S} = \pi^* \circ \mu_{\mathfrak H}$.

Let $[xy]$ denote the (unique) geodesic between points $x, y \in X$.  We first
prove the following Crofton formula.

\begin{proposition}\label {1}                
There is a constant $k > 0$ such that, if $x, y \in X$ then  
$$
\mu_{\mathfrak S} \left \{S \in {\mathfrak S}: S \cap [xy] \not= \emptyset \right \} = k
d(x, y).
$$

\end{proposition}
For the purposes of the proof, we introduce the notation
$$
m(x,y) = \mu_{\mathfrak S} \left \{S \in {\mathfrak S}: S \cap [xy] \not= \emptyset \right \}.
$$

\begin{lemma}\label {2}  If $x, y \in X$ then
 
(a) \quad $\left \{ S \in {\mathfrak S}: S \cap [xy] \not= \emptyset \right \}$ is
compact;

(b) \quad $m(x, y) < \infty$.
\end{lemma}

{\sc Proof:} (a) Let $g, h \in G$.  We claim that the point $gx_0$ lies on the
hyperplane $h S_0$ if and only if $gK \cap h G(S_0) \not= \emptyset$.  The
crucial fact used in the proof is that $G(S_0)$ acts transitively on $S_0$. 
Therefore
\begin{equation*}
\begin{split}
gx_0 \in h S_0 & \Longleftrightarrow h^{-1} g x_0 \in S_0 \\
& \Longleftrightarrow h^{-1} g x_0 = g_0 x_0, \hbox { for some } \ g_0 \in G(S_0) \\
& \Longleftrightarrow g_0^{-1}h^{-1}g \in K, \hbox { for some } g_0 \in G(S_0) \\
& \Longleftrightarrow h^{-1} g = g_0 k_0, \hbox { for some } g_0 \in G(S_0), k_0
\in K \\
& \Longleftrightarrow gK \cap h G(S_0) \not= \emptyset. \\
\end{split}
\end{equation*}
This proves our assertion.  Furthermore, since $[xy]$ is compact, there exists
a compact subset $J \subset G$ such that $Jx_0 = [xy]$, by \cite[p.137, Lemma
1]{nac}.  A hyperplane $hS_0$ meets $[xy]$ if and only if $JK \cap h G(S_0) \not= \emptyset$,
that is $h \in JK G(S_0)$.  The
set of such hyperplanes $h S_0$ is therefore compact, being the image of the
compact set $JK$ under the quotient map $G \rightarrow G/G(S_0)$.

(b) This follows immediately from (a). \qed

\bigskip
{\sc Proof of Proposition \ref {1}}.  We must prove that $m(x, y) = kd(x,y)$
where $k > 0$ is constant.  This is based on the following facts: (a) $X$ is
two-point homogeneous with respect to the action of $G$; (b) the measure
$\mu_{\mathfrak S}$ is $G$-invariant.  It follows that $m(x,y) = m(x^\prime,
y^\prime)$ whenever $d(x,y) = d(x^\prime, y^\prime)$.  Moreover, by considering
a large number of pairwise disjoint geodesic segments of equal length in some
fixed geodesic segment $[ab]$, we see that $m(x,y) \rightarrow 0$ as $d(x,y)
\rightarrow 0$.  In particular $m(x,x) = 0$.   If we divide $[xy]$ into $s$
segments $[x_i x_{i+1}]$ of equal length $(1 \leq i \leq s)$, then $m(x,y) =
sm(x_1, x_2)$.  In this way, if $d(x^\prime, y^\prime) = q d(x,y)$
 where $q$ is rational, then
$m(x^\prime, y^\prime) = q m(x,y)$.  Now $m(x,y) < \infty$, by Lemma 2, and so
by continuity there is a constant $k \geq 0$ such that $m(x,y) = k d(x,y)$ for
all $x, y \in X$.  We must check that $k > 0$.  For this, it is enough to show
that $m(x,y) > 0$ for some $x,y$.  Take $x = (1, 0, 0. \ldots, \epsilon)$, $y =
(1, 0, 0. \ldots, - \epsilon)$.  Then the hyperplane $S_0$ meets $[xy]$
transversally at the interior point $x_0 = (1, 0, \ldots , 0)$ of the geodesic
segment $[xy]$.  There is an open neighbourhood $\widetilde V$ of the identity in $G$
such that $g S_0$ meets $[xy]$ at an interior point for all $g \in \widetilde
V$.  The image $V$ of $\widetilde V$ in $G/G(S_0)$ is an open set which is
contained  in the set of all hyperplanes meeting $[xy]$.  Since $V$ has
positive measure, this implies that $m(x,y) > 0$, as required. \qed 

\begin{remark}
Unlike the case $n=2$ in \cite[Section 3]{s1}, 
\cite[Section 4.4]{gel}, the above proof of the Crofton formula
in hyperbolic space is not constructive. However it gives
a geometric explanation of why the result is true.

\end{remark}

We now re-interpret Proposition 1, replacing $\mu_{\mathfrak S}$ by $\mu_{\mathfrak H}$,
the invariant measure on the space of all half-spaces.  This will allow us to
use the method of \cite{rs} to show that the metric on $X$ is of negative
type.  Given $x \in X$, let $\Sigma_x$ denote the set of half-spaces containing
$x$.

\begin{lemma} \label{3} If $x, y \in X$ then
$$
\mu_{\mathfrak H} (\Sigma_x \bigtriangleup \Sigma_y) = \mu_{\mathfrak S} \left \{ S \in
{\mathfrak S} : S \cap [xy] \not= \emptyset \right \}.
$$
\end{lemma}

{\sc Proof:} The double covering $\pi : {\mathfrak H} \rightarrow {\mathfrak S}$ is given
by $\pi (H) = \partial H$, the boundary of $H$.  Since a hyperplane $S$ is
totally geodesic we have either $\#\left ( [xy] \cap S \right ) \leq 1$ or $[xy]
\subseteq S$.  Therefore
$$
\Sigma_x \bigtriangleup \Sigma_y = \pi^{-1} \left \{ S \in {\mathfrak S} : S \cap [xy] \not=
\emptyset \hbox { and } [xy] \not\subseteq S \right \}.
$$
Now $\mu_{\mathfrak S} \left \{ S \in {\mathfrak S} : [xy] \subseteq S \right \} = 0$.  The result
follows, since $\mu_{\mathfrak S} = \pi^* \circ \mu_{\mathfrak H}$. \qed

\begin{corollary} \label {4}  If $x, y \in X$ then $d(x, y) = k \mu_{\mathfrak H}
(\Sigma_x \bigtriangleup \Sigma_y)$ where $k > 0$ is constant.  Hence $d$ is a
kernel of negative type.
\end{corollary}

{\sc Proof:} The first assertion is immediate from Lemma~\ref{3}.  The second
assertion follows by embedding $X$ into $L^2({\mathfrak H}, \mu_{\mathfrak H})$ via $x
\mapsto v_x$ where $v_x = \chi_x - \chi_{x_{0}}$, and $\chi_x$ is the
characteristic function of $\Sigma_x$.  Then $\mu_{\mathfrak H}(\Sigma_x
\bigtriangleup \Sigma_y) = \| v_x - v_y \|_2^2$ and so $\sqrt d$ is a Hilbert
space distance. (c.f. \cite[Proposition 1.1]{rs}.) \qed

\begin{remark} \label {fh}
Faraut and Harzallah show \cite {far} that the distance function on real or
complex hyperbolic space is a kernel of negative type.  See \cite[Chapitre 6,
Th\'eor\`eme 21]{hav}.  The corresponding result for quaternionic 
hyperbolic space
is false, because the group of isometries $Sp(1, n)$ has Kazhdan's property $(T)$
\cite[Th\'eor\`eme 6.4]{far}. It would be interesting
to have a direct proof of this fact, avoiding the use of property $(T)$.
\end{remark}

\begin{remark} Corollary \ref {4} is stronger than the result of Faraut and
Harzallah.  It asserts that the distance $d$ is a measure definite kernel, in
the sense of \cite{rs}, a concept which is in general strictly stronger than
that of negative type. See the next section.
\end{remark}

\begin {remark}  Since $G(H_0)$ is a closed noncompact subgroup of
$G$, it follows from Moore's ergodicity Theorem \cite[Theorem 2.2.6]{zim} that
any lattice subgroup $\Gamma$ of $G$ acts ergodically on ${\mathfrak H} = G/G(H_0)$. 
Note that the group $\Gamma$ does not have property $(T)$.
Fix an element $x \in X$ then $\mu_{\mathfrak H} (\Sigma_x \bigtriangleup
g\Sigma_x) = \mu_{\mathfrak H} (\Sigma_x \bigtriangleup \Sigma_{gx})$ is a positive
constant multiple of $d(x, gx)$, and hence unbounded.  The action of $\Gamma$ in
 \cite[Theorem 2.1] {rs} may therefore be chosen to be
ergodic.  Whether one can find an ergodic action with similar properties
 for all groups without property $(T)$ is an open question.
\end{remark}

\begin {remark} Consider any Coxeter system $(W,S)$, where $W$ is a Coxeter group
and $S$ is the canonical generating set. 
 Let $\Delta$ denote the associated Coxeter complex.
By a result of J. Tits \cite[Chap. IV, p.41]{b}, the natural distance $d (c, c^\prime)$
between chambers $c$ and $c^\prime$ 
is equal to the number of walls of $\Delta$ separating $c$ and
$c^\prime$. This is precisely the Crofton formula for distance, relative to 
counting measure on the space of walls.

Let ${\mathcal H}$ denote the space of half spaces (``roots'') of $\Delta$ with counting
measure.  If $c$ is a chamber of $\Delta$, let $H_c$ be the set of roots
containing $c$.  Then
$2d (c, c^\prime) = \# (H_c \bigtriangleup H_{c^{\prime}})$.  (See the proof of
\cite [Proposition 6.14]{hav}.)  The argument proceeds as before,
showing that $d$ is a kernel of negative type.

The same type of argument is applied in \cite [Proof of Theorem 3]{bs}
to a polygonal complex $X$,
which is locally finite, simply connected and of type $(4,4)$ or $(6,3)$.
Let ${\mathfrak B}$ denote the set consisting of all barycentres of faces of 
$X$ or of edges of $X$ which are not adjacent to faces.
 An unbounded negative definite kernel $N$
of the above form is defined on ${\mathfrak B}$.  The kernel $N$ is equivalent
to the geodesic distance between elements of ${\mathfrak B}$ and is invariant
under automorphisms of $X$. This is used in \cite{bs} to show that a properly discontinuous group
of automorphisms of $X$ does not have property (T).

\end{remark}

\section{the hypermetric property for $H_{{\RR}}^n$}

A semimetric space $(X,d)$ is said to be {\em hypermetric} if it satisfies the
following property: for each finite subset
$\{x_1,x_2 \dots, x_m\}$ of $X$ and integers $\{t_1,t_2 \dots, t_m\}$
such that $\sum_{i=1}^{m}t_i = 1$, we have
$\sum_{i,j=1}^{m}t_it_jd(x_i,x_j) \le 0$.
The corresponding statement, with $t_1,t_2 \dots, t_m$ real numbers
satisfying $\sum_{i=1}^{m}t_i = 0$, says that $d$ is a kernel negative type.
It is easy to see that if $d$ has the hypermetric property then
$d$ is a kernel negative type \cite[1.2]{a}.

Suppose that $(\mathfrak W, \mu)$ is a measure space and that a semimetric $d$ on a set $X$
is defined  by the formula
 $d(x, y) =  \mu(S_x \bigtriangleup S_y)$, where for each
$x \in X$, $S_x$ is a measurable set.
It was proved in \cite[Theorem 3.1]{ke} that
$d$ is then a hypermetric. The next result is therefore an immediate consequence of
Corollary \ref{4}. 
\begin{corollary} \label {hypermetric}  The geodesic distance $d$ on $H_{{\RR}}^n$
is a hypermetric.
\end{corollary}

Finite hypermetric spaces have been characterized up to isometry in \cite{a},
and there is a detailed exposition of their properties in
\cite {dgl}.
In view of the fact that the space $(H_{{\RR}}^n, d)$ has constant negative curvature,
the next result may seem slightly surprising.

\begin{corollary} \label {spherical}  Let $\{x_1,x_2 \dots, x_m\}$ be any
finite subset of  $H_{{\RR}}^n$, endowed with the metric $\sqrt d$.
 Then $(\{x_1,x_2 \dots, x_m\}, \sqrt d)$
 embeds isometrically in a euclidean sphere in $\RR^p$,
where $p \ge \log_2m$. \end{corollary}

{\sc Proof:} Since $d$ is of negative type (being hypermetric) it follows
from \cite [Proposition 5.14]{hav} that 
$(\{x_1,x_2 \dots, x_m\}, \sqrt d)$ embeds isometrically in $\RR^p$,
for some $p$. The fact that the image is contained in a euclidean
sphere is then a consequence of \cite [Lemme 1.12]{a}. The estimate
$p \ge \log_2m$ is provided by  \cite [Proposition 1.18]{a}.
\qed

\section{complex hyperbolic space}

Complex hyperbolic space $H_{\CC}^n$ is constructed in a manner similar to that
of $H_{\RR}^n$ \cite{chen}.  The bihermitian form $\langle z, w \rangle = -
\overline {z^0} w^0 + \overline {z^1} w^1 + \cdots + \overline {z^n} w^n$ on
${\CC}^n$ defines a set of negative vectors in ${\CC}^{n+1}$,
defined by the condition $\langle z, z \rangle < 0$, which projects to
the subspace $H_{\CC}^n$ of ${\CC}{\PP}^n$.

One might try to mimic the arguments of the preceding section in the complex
case.  The space $H_{\CC}^n$ is two point homogeneous with respect to the
isometry group $U(1, n)$.  Natural analogues of hyperplanes are the
equidistant hypersurfaces \cite[Section 4]{gp}, also known as spinal surfaces in
\cite {mos}.  These are, by definition, subspaces of the form $S =
\left \{ z \in H_{\CC}^n : d(z, z_1) = d(z, z_2) \right \}$, where $z_1, z_2
\in H_{\CC}^n$.  There is a natural invariant measure on the space of
equidistant hypersurfaces, but we cannot prove an analogue of Lemma~\ref{2},
because the group $G(S_0)$ of elements of $U(1, n)$ which leave a given
equidistant hypersurface $S_0$ globally invariant does not act transitively on
$S_0$ \cite [Section 3.2]{mos}.  An additional problem is that equidistant
hypersurfaces are not totally geodesic and a geodesic segment can meet an equidistant
hypersurface more than once, without being contained in the hypersurface. 
However, even in the complex case, it is possible to prove a Crofton formula if
we consider horospheres instead of equidistant hypersurfaces.

Horospheres in complex hyperbolic space are described geometrically in \cite
[Section 1]{gp}.  From a group-theoretic point of view they may be described as
follows \cite[Chapter II.1]{hel}.  Let $G = U(1, n)$, so that $H_{\CC}^n = G/K$ where
$K = U(1) \times U(n)$.  Let $G = KAN$ be the Iwasawa decomposition and $M$ the
centralizer of $A$ in $K$.  Thus $MAN$ is a minimal parabolic subgroup.    Let
$x_0$ be the origin $K$ in $G/K$ and $\xi_0 = Nx_0$  Then $\xi_0$ is a
horosphere and the subgroup of $G$ which maps $\xi_0$ to itself equals $MN$. 
The homogeneous space $G/MN$ is the space of horocycles.  Note that $N$ is
isomorphic to the complex Heisenberg group \cite [Chapter 2, p.215]{hel}.  The
groups $G$ and $MN$ are therefore unimodular \cite[Chapter X, Proposition
1.4]{h1}, and there is a $G$-invariant measure on $G/MN$.  Since $MN$ acts
transitively on $\xi_0$, the set of horospheres which meet a geodesic segment
$[xy]$ is compact, as in Lemma~\ref{2}.  Consider the expression
$m(x,y) = \int_{G/MN} n(h) d\mu (h)$, where $n(h)$ is the number of times a
horocycle $h \in G/MN$ meets the geodesic segment $[xy]$. Since $H_{\CC}^n$ is
two-point homogeneous, and $\mu$ is $G$-invariant, $m(x,y)$ depends only on
$d(x,y)$.  The same argument as in the proof of Proposition~\ref {1}
establishes the following Crofton formula.

\begin {proposition} \label {6}  There is a constant $k > 0$ such that, for all
$x, y \in H_{\CC}^n$,
$$
\int_{G/MN} n[x,y] (h) d \mu (h) = k d(x,y)
$$
where $n[x,y] (h)$ is the number of times that $h$ meets $[xy]$.
\end{proposition}

\begin{remark}  Because it is possible to have $1 < n[x,y] (h) < \infty$, it is
not clear how to use this result to prove that $d(x,y)$ is of negative type, as
we did for the real case in Corollary \ref{4}. In fact, the example
in the next section indicates that a Crofton formula alone is not enough.
Note that a
Crofton formula involving horospheres is clearly also valid in real hyperbolic space.
\end{remark}

\section{projective space}

The purpose of this section is twofold. Firstly we give an example to show that
the existence of a Crofton formula for a distance function does not in general imply
that the distance is of negative type. Secondly, this example completes the 
classification of riemannian symmetric spaces of rank one for which the geodesic
distance is of negative type.

The projective spaces $P_{\FF}^n$ of dimension $n \ge 2$ over ${\RR}$, ${\CC}$, ${\HH}$ and (in
the case $n = 2$) the octonions ${\mathbb O}$, are compact two-point homogeneous spaces.

\begin{proposition} The geodesic distance $d$ on $P_{\FF}^n$ is not of
negative type.
\end{proposition}

The projective spaces all contain $P_{\RR}^2$ as a geodesic subspace,
 and so it is enough to consider the case
of $P_{\RR}^2$.
As usual, the same notation is used for a vector $x \in {\RR}^{3}$ and its
equivalence class in $P_{\RR}^2$.  The geodesic distance between points $x,
y$ in  $P_{\RR}^2$ is given by $\cos d(x, y) = \frac{\left |
( x, y ) \right |}{ \left ( ( x, x ) ( y, y) \right )^{\frac{1}{2}}}$, where $(x, y)$ is the
usual inner product on ${\RR}^{3}$.

\begin{proposition} The geodesic distance $d$ on $P_{\RR}^2$ is not of
negative type.
\end{proposition}

{\sc Proof:} We show that the
metric $d$ is not of negative type by exhibiting points
$x_1,x_2 \dots,x_6 \in P_{\RR}^2$ and  real numbers $t_1,t_2 \dots, t_6$
satisfying $\sum_{i=1}^{6}t_i = 0$ such that $\sum_{i,j=1}^{6}t_it_jd(x_i,x_j) > 0$.
In fact we choose $(t_1,t_2 \dots, t_6) = (1,1,1,-1,-1,-1)$ and
for notational convenience put
$(x_1,x_2,x_3,x_4,x_5,x_6) = (p_1,p_2,p_3,q_1,q_2,q_3)$. It is enough to choose these
points so that
\begin{equation}\label{inequality}
\sum d(p_i,p_j) + \sum d(q_i,q_j) > \sum d(p_i,q_j).
\end{equation}
We make the following choices.
Let $p_1~=~(1,0,1)$, $p_2~=~(1,0,-1)$, $p_3~=~(0,1,0)$, and
 $q_1~=~(0,1,1)$, $q_2~=~(0,1,-1)$, $q_3~=~(1,0,0)$.
Then $d(p_i,p_j) = \pi/2$, $d(q_i,q_j) = \pi/2$, for $i,j \in \{1,2,3\}$.
Moreover

\begin{equation*} d(p_i,q_j) = 
\begin{cases}
\pi/3 & \text{if  both $i,j \in \{1,2\}$},\\
\pi/4  & \text{if exactly one of  $i,j$ equals $3$},\\ 
\pi/2  & \text{if $i=3$ and $j=3$}.
\end{cases}
\end{equation*}

We therefore obtain
$\sum d(p_i,p_j) + \sum d(q_i,q_j) = 3\pi/2 + 3\pi/2 = 3\pi$ and
$\sum d(p_i,q_j) = 4\pi/3 +4\pi/4 + \pi/2 = 17\pi/6$, which proves 
the inequality (\ref{inequality}).
\qed

\bigskip
{\sc Example.}\ There is a Crofton formula for geodesic 
distance $d$ on the space $P_{\RR}^n$, but as we have just seen,
this distance is not of negative
type. The proof of a Crofton formula for $d(x,y)$ in terms of totally geodesic
hypersurfaces which meet the segment $[xy]$ is exactly the same as for real
hyperbolic space in Section \ref{rhyperbolic}, except that by
compactness there is no need for an analogue of Lemma \ref {2}.
The reason that the analogue of Corollary \ref{4} fails is that a
totally geodesic hypersurface does not separate $P_{\RR}^n$
into two parts. 

\begin{remark} Two point homogeneous riemannian manifolds have been
classified completely \cite[Chapter IX \S5]{h1}. They are the euclidean
spaces, the circle, and the symmetric spaces of rank one. The symmetric spaces of 
rank one for which the geodesic distance is of negative type are the spheres
and the real or complex hyperbolic spaces. This follows from the results
of \cite{far}, together with our result for projective spaces. Note that it
is easy to give a new proof that the geodesic distance on a sphere is of negative type
along the lines of Corollary \ref{4}, using a measure on the space of
half-spaces. 
\end{remark}

\begin{remark}
There is a naturally occurring invariant metric of negative type on
the  projective space $P_{\FF}^n$ \cite{hi}. For there is an embedding of
$P_{\FF}^n$ as the set of primitive idempotents in the formally real Jordan algebra ${\mathfrak A}$
of hermitian $n \times n$ matrices over $\FF$.
The trace form defines a euclidean distance $d_E$ on ${\mathfrak A}$
which induces a metric of negative type on $P_{\FF}^n$.
It follows from \cite {hi} that $P_{\FF}^n$ is
also two point homogeneous relative to this metric. 
\end{remark}

\bigskip
\noindent {\bf ADDENDUM}

In Remark 2.6 we asked for a direct proof
of the fact that the distance
function on quaternionic hyperbolic space is a kernel of negative type.
This is in fact possible, by considering the following 24 points.
For $\sigma \in \{+1,-1\}$ and $\epsilon \in \{\pm i, \pm j,\pm k\}$, let
$x^{\sigma}_{\epsilon}=(3,2\sigma + 2\epsilon,0)$ and
$y^{\sigma}_{\epsilon}=(3,0,2\sigma + 2\epsilon)$.
Direct calculation shows that
$$\sum d(x^{\sigma}_{\epsilon},x^{\rho}_{\delta})
+ \sum d(y^{\sigma}_{\epsilon},y^{\rho}_{\delta})
- \sum d(x^{\sigma}_{\epsilon},y^{\rho}_{\delta})
= 417.03 - 415.77 >0.$$
This shows that the condition for negative type fails with the number $+1$
being associated to the $x^{\sigma}_{\epsilon}$ and $-1$ to the
$y^{\sigma}_{\epsilon}$.

\end{document}